\documentclass[a4paper,11pt]{article} \usepackage{amsmath,amssymb,amsthm}
\usepackage[english]{babel}

\newtheorem{theo}{Theorem} 
 \newtheorem{prop}[theo]{Proposition}

\newcommand{\Na}{\mathbb N}                   

\newcommand{\Za}{\mathbb Z}                   

\newcommand{\Ra}{\mathbb R}                   

\newcommand{\scal}[1]{\langle #1 \rangle}

\newcommand{\finpreuve}{\hfill $\Box$}

\newcommand{\name}{$\underline{\qquad \qquad}$} 

\begin{document}

\author{  Jean-Marc Bouclet\footnote{Jean-Marc.Bouclet@math.univ-toulouse.fr} and Yannick Sire\footnote{sire@cmi.univ-mrs.fr}}
\title{{\bf \sc Refined Sobolev inequalities on  manifolds with ends}}

\maketitle

\begin{abstract} By considering a suitable Besov type norm, we obtain refined Sobolev inequalities on a family of Riemannian manifolds with (possibly exponentially large) ends. The interest is twofold: on one hand, these  inequalities are stable by multiplication by rapidly oscillating functions, much as the original ones \cite{GMO}, and on the other hand our Besov space is stable by spectral localization associated to the Laplace-Beltrami operator (while $ L^p $ spaces, with $ p \ne 2 $, are in general not preserved by such localizations on manifolds with exponentially large ends). We also prove an abstract version of refined Sobolev inequalities for any selfadjoint operator on a measure space (Proposition \ref{general}). 
  
  
  
\end{abstract}






\setcounter{equation}{0}
 For functions $ u$ on $ \Ra^n $, $ n \geq 1 $, and $ p \geq 2 $ real, the homogeneous and inhomogeneous Sobolev estimates
\begin{eqnarray}
 || u ||_{L^{p} (\Ra^n)} & \leq & C || (-\Delta)^{s_p/2} u ||_{L^2 (\Ra^ n)} =: || u ||_{\dot{H}^{s_p}} \nonumber \\
 & \leq &  C || (1-\Delta)^{s_p /2} u ||_{L^2(\Ra^n)} =: || u ||_{H^{s_p}} , \qquad s_p  = \frac{n}{2} - \frac{n}{p}  , \label{pournotationsp}
\end{eqnarray}
 are a very well known tool to control $ || u ||_{L^{p}} $  when $u$ is sufficiently smooth.
The homogeneous version  is sharp with respect to all scalings  and the inhomogeneous one is sharp with respect to high frequency scalings ({\it i.e.} $ u (x) \mapsto u (\lambda x) $ with $ \lambda \geq 1  $). The drawback of these estimates is to behave
  badly under the multiplication by  characters. The usual counterexample (see {\it e.g.} \cite[subsection 1.3.2]{BCD}) is to consider
\begin{eqnarray}
  u_{\varepsilon} (x) = e^{\frac{i}{\varepsilon} x \cdot \eta } \phi (x) , \label{counterexampleRn}
\end{eqnarray}
with $ \phi \ne 0 $ in the Schwartz space $ {\mathcal S} (\Ra^n) $, $ \eta \ne 0 $ in $ \Ra^n  $, and to observe that
$$ || u_{\varepsilon} ||_{L^{p}} = || \phi ||_{L^{p}}, \qquad || u_{\varepsilon} ||_{\dot{H}^{s_p}} \sim || u_{\varepsilon} ||_{H^{s_p}} \sim \varepsilon^{-s_p} |\eta|^{s_p} || \phi ||_{L^2} ,  $$
for which the Sobolev estimates alone  provide only the very bad estimate $ || u_{\varepsilon} ||_{L^{p}} \lesssim \varepsilon^{-s_p} $. Notice that  considering the homogeneous or inhomogenous Sobolev norm is irrelevant here. A nice well known substitute to the Sobolev estimates is given by the following so called refined Sobolev inequalities (introduced in \cite{GMO}, see also \cite{BCD} for references)
\begin{eqnarray}
 || u ||_{L^p} & \leq & C_p || u ||_{\dot{H}^{s_p}}^{\frac{2}{p}} || u ||_{\dot{B}^{s_p-\frac{n}{2}}}^{1-\frac{2}{p}} \label{substitut1} \\
 & \leq & C_p || u ||_{H^{s_p}}^{\frac{2}{p}} || u ||_{\dot{B}^{s_p-\frac{n}{2}}}^{1-\frac{2}{p}} , \label{substitut0}
\end{eqnarray} 
where, for $ \sigma > 0 $, $ \dot{B}^{-\sigma} = \dot{B}^{-\sigma}(\Ra^n) $ is the homogeneous Besov space with norm
\begin{eqnarray}
 || u ||_{\dot{B}^{-\sigma}} = \sup_{ \lambda > 0} \lambda^{-\sigma} || \theta \big(\lambda^{-1}(-\Delta)^{1/2} \big) u ||_{L^{\infty}(\Ra^n)} ,\label{Besovnorm}
\end{eqnarray} 
with $ \theta \in C_0^{\infty}(\Ra^n) $  a fixed bump function such that $ \theta \equiv 1 $ near $ 0 $. For functions of the form $ u_{\varepsilon} $, one can easily see 
that, for any $ \sigma \in (0,n] $, $
 || u_{\varepsilon} ||_{\dot{B}^{-\sigma}} \leq C \varepsilon^{\sigma} $
hence if we choose $ \sigma = \frac{n}{2} - s_p $, (\ref{substitut0}) yields the sharp\footnote{with respect to $ \varepsilon \ll 1 $} bound 
\begin{eqnarray}
 || u_{\varepsilon} ||_{L^p} \lesssim \varepsilon^{-\frac{2 s_p}{p}} \varepsilon^{\left( \frac{n}{2} - s_p \right)  \left( 1 - \frac{2}{p} \right)} \lesssim 1 . \label{familyRn}
\end{eqnarray}
We recall that refined Sobolev inequalities are also useful to  study the lack of compactness of Sobolev embeddings as in \cite{Gerard}. For such applications, we also refer to the recent paper \cite{PP} (and the references therein) where improved Sobolev inequalities (in the Morrey scale) are used to restore the compactness of maximizing
sequences for the Sobolev inequality up to dilation and translation. The lack of compactness in Sobolev embeddings is in turn related to the profile decomposition  of solutions to some PDE; we refer for instance to \cite{Laurent} where refined Sobolev inequalities are explicitly used in this context. See also the book
\cite{FT}.


The purpose of this  note is to give a robust analogue of (\ref{substitut0}) on a class of Riemannian manifolds with ends. By robust we mean in particular that the related Besov norms should be stable by spectral localization, which in practice can be interesting to combine those estimates with a Littlewood-Paley decomposition or focus only on the high frequency part of the function ({\it i.e.} $  \rho (h^2 \Delta_g) v $ with $ h \ll 1 $ and $ \rho \in C^{\infty} (\Ra) $ equal to $0 $ near $ 0 $). Unfortunately, in many natural cases as manifolds with exponentially growing ends ({\it e.g.} the hyperbolic space), $ L^p $ spaces (in particular $ L^{\infty} $) are in general not preserved by  $ \chi ( \Delta_g ) $ when $ \chi $ belongs to $ C_0^{\infty} $ so the definition (\ref{Besovnorm}) (with $ \Delta_g $ instead of $ \Delta $) does not seem to be adapted to such a purpose. We will introduce below alternative norms for which  refined Sobolev inequalities still hold and which
  are sharp enough  to handle rapidly oscillating functions as in (\ref{counterexampleRn}).

Before turning to this specific question we consider first an abstract general extension of (\ref{substitut0}) in Proposition \ref{general}. We will use it in the special case of manifolds with ends but, since it is very simple and (as we feel) of independent interest, we prefer to isolate it first.
The generalization of refined Sobolev inequalities has already been considered in \cite{Ledoux} where Ledoux proved  that, up to the replacement of $ \theta (\lambda^{-1} (-\Delta)^{1/2}) $ (in the definition of the Besov norm (\ref{Besovnorm})) by families of operators satisfying pseudo-Poincar\'e inequalities (such as the heat semigroup or local averages on balls), refined Sobolev inequalities still hold on various domains carrying a natural gradient or second order operator (Riemannian manifolds with Ricci curvature bounded from below, Cayley graphs).  Furthermore, his approach allows him to study the dependence of the constants on the dimension of the domain. 
Our first remark, summarized in Proposition \ref{general} below, is in this spirit and is basically the simple observation that the proof of \cite[Theorem 1.43]{BCD} is completely decorrelated to any Laplacian type operator. In other words,  refined $ L^2 $ Sobolev inequalities rest on essentially nothing and hold with some kind of 'universal' bound (depending only on $p$ and the cutoff $ \theta $ involved in the relevant Besov (type) norm). Here is the precise statement.

 Let $ \big( {\mathcal M} , \mu \big) $ be a measure space. For $ p \in [1,\infty] $, we set  $ L^p = L^p ({\mathcal M},d \mu) $. Let $ (A,\mbox{Dom}(A)) $ be a nonnegative selfadjoint operator on $ L^2 = L^2 ({\mathcal M},d\mu) $ and, for $ s > 0 $,
  define
\begin{eqnarray}
  || u ||_{\dot{H}^s_A} := | | A^{s} u | |_{L^2} , \label{homogeneous}
\end{eqnarray}
 on $ \mbox{Dom} (A^{s})  $.
We next let $ \theta \in C_0^{0}(\Ra) $ be such that $ \theta \equiv 1 $ near $ 0 $ and set, for $ \sigma \geq 0 $,
\begin{eqnarray}
  || u ||_{\dot{B}_A^{-\sigma}} := \sup_{k \in \Za} 2^{-k \sigma}|| \theta (2^{-k} A) u ||_{L^{\infty}} , \label{definitionBesov}
\end{eqnarray}  
on  the subspace of elements $ u \in L^2 $ such that $ || u ||_{\dot{B}_A^{-\sigma}} $ is finite.

\begin{prop} \label{general} Assume that $ \theta \equiv 1 $ on $ [-c,c] $. Then for all  $ s  > 0 $, $  \sigma > 0 $ and 
\begin{eqnarray}
 p := \frac{2(\sigma + s)}{\sigma}  , \label{pSobolev}
\end{eqnarray}
we have 
$$ || u ||_{L^p} \leq C  || u ||^{1-\frac{2}{p}}_{\dot{B}_A^{-\sigma}} || u ||^{\frac{2}{p}}_{\dot{H}_A^s} , $$
 for all $ u \in \emph{Dom} (A^s)  $ such that $ || u ||_{\dot{B}_A^{-\sigma}} $ is finite. The constant is explicitly given by
 $$ C = 2 || 1 - \theta ||_{L^{\infty}}^{\frac{2}{p}} \left( \frac{p}{p-2} \right)^{1/p} \left(\frac{2}{c} \right)^{\sigma \frac{p-2}{p}} . $$
\end{prop}

Let us point out that the purely spectral nature of  refined $ L^2 $ Sobolev inequalities was already observed in \cite{Ledoux}, but in a slightly different and more complicated form. Here we adapt the proof of \cite[Theorem 1.43]{BCD} so don't claim much originality. Our point is only to emphasize
 that no assumption on $A$ (but its selfadjointness) is required. In particular,  this result has nothing to do with any Laplacian type operator (continuous or discrete), as one can see for instance from (\ref{pSobolev}) which does not involve any dimension. Of course, in practice, the finiteness of $ || u ||_{\dot{B}^{-\sigma}_A} $ can be checked only on a case by case basis, using  a non trivial knowledge of functions of $A$, but Proposition \ref{general} shifts  the proof of refined Sobolev inequalities to suitable estimates on 
$ || \theta (2^{-k} A) u ||_{L^{\infty}} $  in a simple and explicit fashion.

We finally point out that, if $( {\mathcal M} , g ) $ is a compact Riemannian manifold with measure $ \mu = {\rm vol}_g $, and if we set $ A := (-\Delta_g)^{1/2} $ then $ || 1 ||_{\dot{H}_A^s } = 0 $ for all $s> 0$  while $ \theta (2^{-k} A) 1 = 1 $ for all $k$ so that 
$$ || 1 ||_{\dot{B}_A^{-\sigma}} = + \infty , $$
 due to  those $ k $ going to $ - \infty $. Therefore, in this case, Proposition \ref{general} gives no contradiction with the fact that $ || 1 ||_{L^p} = \mu ({\mathcal M})^{1/p} > 0 $ and $ || 1 ||_{\dot{H}^s_A} = 0 $.
 
\bigskip

\noindent {\it Proof.} One starts with the usual formula
\begin{eqnarray}
 || u ||_{L^p}^p = p \int_0^{\infty} \lambda^{p-1} \mu \big( \{ |u| > \lambda \} \big) d \lambda . \label{usualformula}
\end{eqnarray}
Assuming that  $ ||u||_{\dot{B}_A^{-\sigma}}> 0 $ as we may (otherwise $ u = 0 $ a.e. and the result is trivial),  we choose the unique integer $ k_{\lambda}  $ satisfying
\begin{eqnarray}
 2^{\sigma k_\lambda} || u ||_{\dot{B}_A^{-\sigma}} \leq \frac{\lambda}{2} < 2^{\sigma (k_\lambda+1)} || u ||_{\dot{B}_A^{-\sigma}}  . 
 \label{conditiosurk2}
\end{eqnarray}
%
We next introduce the decomposition of $u = v_{\lambda} + w_{\lambda}$ with
$$ v_{\lambda} := \theta (2^{-k_{\lambda}}A)  u,\qquad w_{\lambda} = \big( 1 - \theta (2^{-k_{\lambda}}A) \big) u . $$
By definition of $ || u ||_{\dot{B}_A^{-\sigma}} $, we have
 $ || v_{\lambda} ||_{L^{\infty}} \leq 2^{\sigma k_\lambda} || u ||_{\dot{B}_A^{-\sigma}} $ hence
$$ \mu \big( \{ |v_{\lambda}| > \lambda/2 \} \big) = 0 . $$
On the other hand, the Markov-Tchebichev inequality yields
$$ \mu \big( \{ |w_{\lambda}| > \lambda/2 \} \big) \leq \frac{4}{\lambda^2} || w_{\lambda} ||_{L^2}^2 , $$
and, if we let $ (E_t)_{t \in \Ra} $ be the resolution of identity associated to $ A $,   the Spectral Theorem yields
$$ || w_{\lambda} ||_{L^2}^2 = \int_{\Ra}  \big| 1 - \theta (2^{-k_{\lambda}}t) \big|^2 d (u, E_t u) . $$
If $ \theta \equiv 1 $ on $ [-c,c] $, we  obtain
$$ || w_{\lambda} ||_{L^2}^2 \leq \sup |1-\theta|^2 \int_{|t|> c 2^{k_{\lambda}}}  d (u, E_t u)  \leq  
\sup |1-\theta|^2 \int_{|t|> \frac{c}{2} \left( \frac{\lambda}{2} \right)^{\frac{1}{\sigma}} || u ||^{-\frac{1}{\sigma}}_{\dot{B}_A^{-\sigma}}}  d (u, E_t u) $$
the second inequality following from (\ref{conditiosurk2}). Coming back to (\ref{usualformula}), we obtain
\begin{eqnarray*}
  || u ||_{L^p}^p & \leq & 4 p \sup |1-\theta|^2 \int_0^{\infty} \lambda^{p-3} \left( \int_{|t|> \frac{c}{2} \left( \frac{\lambda}{2} \right)^{\frac{1}{\sigma}} || u ||^{-\frac{1}{\sigma}}_{\dot{B}_A^{-\sigma}}}  d (u, E_t u)  \right) d \lambda \\
  & \leq & 4 p \sup |1-\theta|^2 \int_{\Ra} \left( \int_0^{ 2(\frac{2}{c})^{\sigma} || u ||_{\dot{B}_A^{-\sigma}} |t|^{\sigma} } \lambda^{p-3}  d \lambda\right)
  d (u, E_t u)
\end{eqnarray*}
where the integral in $ \lambda $ equals $ \frac{2^{p-2}}{p-2} \left( \frac{2}{c} \right)^{\sigma(p-2)} || u ||_{\dot{B}_A^{-\sigma}}^{p-2} |t|^{\sigma(p-2)}  $.
On the other hand, using (\ref{pSobolev}) we have
$$ \int |t|^{\sigma(p-2)} d (u,E_t u) = \int |t|^{2s} d (u,E_t u) = \big| \big| A^s u  \big| \big|^2_{L^2} , $$
and the result follows. \finpreuve

\bigskip

We now consider the case  of manifolds with ends. We will work with the same class of manifolds  as in \cite{Bo1,Bo2}. We briefly recall their definition. We assume that $  ({\mathcal M},g)  $ is a smooth Riemannian manifold without boundary which is diffeomorphic outside a compact subset to a product $ (R,\infty) \times \Gamma $, for some $ R > 0 $ and some compact manifold without boundary $ \Gamma $, such that the metric $g$  reads on  coordinate patches of the form $ (R,\infty) \times U $, with $ U $ a coordinate patch on $ \Gamma $,
$$ g = G (r,y , dr , w (r)^{-1} d y ) , $$
with $ G (r,\sigma, v_1 , v^{\prime}) $ a uniformly elliptic polynomial of degree $2$ in $ (v_1 , v^{\prime})$, with coefficients bounded (as well as their derivatives) with respect to $ (r, y) $. The smooth function $w$, which determines the type of ends, satisfies the conditions
$$ 0 < w (r) \leq C, \qquad  w (r) / w (r^{\prime}) \leq C \ \ \mbox{for} \ \ |r^{\prime}-r| \leq 1, \qquad |\partial_r^k w (r) | \leq C_k w (r) , $$
for all $ r,r^{\prime} > R $. For instance, these assumptions are satisfied by the hyperbolic space $ {\mathbb H}^n $  with  $ w (r) = e^{-r} $ and by $ \Ra^n $ with $ w (r) = r^{-1} $.

In the sequel, we let
$$ A = (- \Delta_g )^{1/2} , \qquad d \mu = d {\rm vol}_g ,  $$
and consider a function $ \theta \in C_0^{\infty} (\Ra) $ such that $ \theta \equiv 1 $ near $ 0 $. For a given $s \in \Ra $, we set
$$ || u ||_{H^s} = || (1-\Delta_g)^{s/2} u ||_{L^2} = || (1+A^2)^{s/2} u ||_{L^2}, $$
and for a given $ \sigma \geq 0 $, we set
\begin{eqnarray}
 || u ||_{\widetilde{B}^{-\sigma}} := \max \left(  \sup_{k \geq 0} 2^{- k \sigma} || \theta (2^{-k}A) u ||_{L^{\infty}}, || u ||_{H^{- \sigma}}  \right) . \label{definitionnorme}
\end{eqnarray} 
Note the difference with (\ref{Besovnorm}): we replace the contribution of low frequencies, {\it i.e.} $ k \leq 0 $ (or $ 0 < \lambda < 1 $), by the Sobolev norm $ || u ||_{H^{-\sigma}} $.
We define $ \widetilde{B}^{-\sigma} $ as the subspace of all $u$ in $ L^2 $ such that $ || u ||_{\widetilde{B}^{-\sigma}} < \infty $.
\begin{prop} \label{pourfinitude} For all $ \chi \in C_0^{\infty}(\Ra) $, 
\begin{enumerate}
\item{the space $ \widetilde{B}^{-\sigma} $ is stable by $ \chi (\Delta_g) $,}
\item{There exists a constant $ C > 0 $ such that, for all $ u \in \widetilde{B}^{-\sigma} $ and all $ j \in \Na $,
$$ || \chi (2^{-2j}\Delta_g) u ||_{\widetilde{B}^{-\sigma}} \leq C || u ||_{\widetilde{B}^{-\sigma}} , $$
and
\begin{eqnarray}
 || \chi (2^{-2j}\Delta_g) u ||_{L^{\infty}} \leq C 2^{j \sigma} || u ||_{\widetilde{B}^{-\sigma}} .  \label{localization}
\end{eqnarray} }
\item{If we replace $ \theta $ in (\ref{definitionnorme}) by another $\theta_1 \in C_0^{\infty}(\Ra) $ such that $ \theta_1 \equiv 1 $ near $ 0 $, we obtain an equivalent norm and the same space $ \widetilde{B}^{-\sigma} $.}
\end{enumerate}
\end{prop}

The items 1 and 2 show that the 'Besov norm' $ || \cdot ||_{\widetilde{B}^{-\sigma}} $ is stable by spectral localization and has the natural behaviour with respect to the semiclassical localizations  $ \chi (2^{-2j} \Delta_g) $.  

\bigskip

\noindent {\it Proof.} We prove the items 1 and 2 simultaneously. If $ u$ belongs to $            \widetilde{B}^{-\sigma} $ then  $ \chi (2^{-2j}\Delta_g) u $ belongs to $ L^2 $. By the Spectral Theorem, we have obviously
$$ || \chi (2^{-2j}\Delta_g) u ||_{H^{-\sigma}} \leq \big( \sup |\chi| \big) \big| \big| u \big| \big|_{H^{-\sigma}}  .$$
On the other hand, it follows from Theorem 2.1 of \cite{Bo1} that, for all $ v \in L^2 \cap L^{\infty} $,
\begin{eqnarray}
 || \chi (2^{-2j} \Delta_g) v ||_{L^{\infty}} \leq C \big( || v ||_{L^{\infty}} + ||v||_{H^{-\sigma}} \big) , \label{finitude}
\end{eqnarray}
and therefore, by taking $ v  =\theta (2^{-k}A) u $ with  $ k \geq 0$, we get 
$$ || \theta (2^{-k}A) \chi (2^{-2j}\Delta_g) u ||_{L^{\infty}} = || \chi (2^{-2j}\Delta_g) \theta (2^{-k}A) u ||_{L^{\infty}} \leq C \big(  || \theta (2^{-k}A) u ||_{L^{\infty}} + ||u||_{H^{-\sigma}}\big)  .$$
Since $ \sigma \geq 0 $, this implies that
$$ \sup 2^{-k\sigma} || \theta (2^{-k}A) \chi (2^{-2j}\Delta_g) u ||_{L^{\infty}} \leq C || u ||_{\widetilde{B}^{-\sigma}} , $$
which proves the first item and the first estimate of the item 2. Let us now prove (\ref{localization}). Since $ \theta \equiv 1 $ near $ 0 $, there exists $ M > 0 $ such that 
$$ \chi (2^{-2j} \lambda^2) \theta (2^{-k} \lambda) = \chi (2^{-2j} \lambda^2), \qquad k - j \geq  M .  $$
Therefore, if we take $ k = j + M $,
\begin{eqnarray*}
 || \chi (2^{-2j}\Delta_g)  u ||_{L^{\infty}} & = & || \chi (2^{-2j}\Delta_g) \theta (2^{-j-M}A) u ||_{L^{\infty}} \\
 & \leq &  C \left(|| \theta (2^{-j-M}A) u ||_{L^{\infty}} + || \theta (2^{-j-M}A)  u ||_{H^{-\sigma}} \right) \\
 & \leq & C 2^{(j+M )\sigma} || u ||_{\widetilde{B}^{-\sigma}} ,
\end{eqnarray*} 
the second line following from (\ref{finitude}) and the third one from the definition (\ref{definitionnorme}).
This proves (\ref{localization}). The item 3 is a simple consequence of (\ref{localization}) with $ \chi (\lambda) = \theta_1 (|\lambda|^{1/2}) $, since if we denote by $ || \cdot ||_{\widetilde{B}_1^{-\sigma}} $ the norm relative to $ \theta_1 $, we have
$$ || u ||_{\widetilde{B}_1^{-\sigma}}  \leq C || u ||_{\widetilde{B}^{-\sigma}}  $$
by (\ref{localization}). Since $ \theta $ and $ \theta_1 $ play symmetric roles, we obtain the expected equivalence of norms.  \finpreuve

\bigskip
In the following theorem, $n$ denotes the dimension of $ {\mathcal M} $.
\begin{theo}[Refined Sobolev inequalities] \label{maintheorem} Let $  s \in (0,n/2) $ and set
$$ p = \frac{2n}{n-2s}, \qquad \sigma = \frac{n}{2} - s > 0 . $$
Then, there exists $ C_p $ such that,
$$ || u ||_{L^p} \leq C_p || u ||_{\widetilde{B}^{-\sigma}}^{1-\frac{2}{p}} || u ||_{H^s}^{\frac{2}{p}} $$
 for all $ u \in L^2 = L^2 ({\mathcal M},d {\rm vol}_g) $ such that the right hand side if finite.
\end{theo}

\bigskip

\noindent {\it Proof.} Let $ \chi \in C_0^{\infty} $ be equal to $ 1 $ near $ 0 $. Using the rough Sobolev embedding (see \cite{Bo2}, more precisely by combining (1.36) in Theorem 1.5 and Lemma 2.4  therein)
$$ || v ||_{L^p} \leq C || v ||_{H^{\frac{n}{2}+1}} $$
and the Spectral Theorem, we have
\begin{eqnarray}
 || \chi (A) u ||_{L^p} \leq C || \chi (A) u ||_{H^{\frac{n}{2}+1}} \leq C^{\prime}  || \chi (A) u ||_{H^{-\sigma}} \leq C^{\prime \prime} || u ||_{H^{-\sigma}}^{1-\frac{2}{p}}
||u||_{H^s}^{\frac{2}{p}} .   \label{ref}
\end{eqnarray}
On the other hand, we may choose $ \chi $ such that $ \theta (2^{-k} A) (1-\chi)(A) = 0 $ if $ k < 0 $, hence using the definition (\ref{definitionBesov}),
$$ || (1-\chi)(A) u ||_{\dot{B}_A^{-\sigma}} = \sup_{k \geq 0} 2^{-k\sigma}|| \theta (2^{-k}A) (1 - \chi)(A) u ||_{L^{\infty}} \leq ||(1-\chi)(A) u ||_{\widetilde{B}^{-\sigma}} $$
whose right hand side is finite by Proposition \ref{pourfinitude} and the fact that $ || u ||_{\widetilde{B}^{-\sigma}} $ is finite. Here we use additionally  that $ \chi (A) = \widetilde{\chi} (\Delta_g) $ with $ \widetilde{\chi}(\lambda) = \chi (|\lambda|^{1/2}) $ which is smooth since $ \chi \equiv 1 $ near $0$. Using Proposition \ref{general}, we thus obtain
\begin{eqnarray*}
 ||(1-\chi)(A)u||_{L^p} & \leq & C_p  || (1-\chi)(A) u ||_{\dot{B}^{-\sigma}_A}^{1 - \frac{2}{p}} || (1-\chi)(A) u ||_{H^s_A}^{\frac{2}{p}} \\
 & \leq & C  ||  u ||_{\widetilde{B}^{-\sigma}}^{1 - \frac{2}{p}} || u ||_{H^s}^{\frac{2}{p}}
\end{eqnarray*}
which, combined with (\ref{ref}), yields the result.
 \finpreuve

\bigskip

We next check that the refined Sobolev inequalities of Theorem \ref{maintheorem} are sufficiently well designed to handle the same kind of counterexamples as (\ref{counterexampleRn}).
Let $U$ be an open coordinate chart of the angular manifold $ \Gamma $ and $ \gamma \in C_0^{\infty}(U) $. We define
$$ u (r,m) = w(r)^{\frac{n-1}{2}}  \psi (r ) \gamma (m), \qquad r \geq R_1, \ m \in \Gamma,  $$
with $ \psi \in {\mathcal S}(\Ra) $ such that $ \mbox{supp}(\psi) \subset [R_1 , \infty ) $, and then
\begin{eqnarray}
 u_{\varepsilon} := e^{i \frac{r}{\varepsilon}} u , \label{family}
\end{eqnarray} 
which is globally defined on $ {\mathcal M} $.
The normalizing factor $ w(r)^{(n-1)/2} $ is here to guarantee that $u$ belongs at least to $ L^2 $ (if $w^{-1}$ grows exponentially, being in the Schwartz space in $r$ is not sufficient to be in $ L^2 $).
Actually $ u $ belongs to all $ L^p $ spaces with $ p \geq 2 $ since 
$$ || u ||_{L^p} \leq C || \gamma ||_{L^p(\Gamma)} \left( \int |\psi(r)|^p w(r)^{(n-1) \left(\frac{p}{2} - 1 \right)} dr \right)^{1/p} , $$
where the right hand side is finite since $w$ is bounded. We also note that the separation of variables in the definition of $ u $ is  irrelevant, and is only for notational convenience in the proof below.

Using the pseudodifferential calculus of \cite{Bo2} (see also \cite[Theorem 2.1]{Bo1}), we will prove the following result.

\begin{prop} \label{contreexemplevariete} Let $ \sigma \geq 0 $ and $ s \geq 0 $. There exists $ C> 0 $ such that, for all $ 0 < \varepsilon < 1 $,
\begin{eqnarray}
|| u_{\varepsilon} ||_{H^{-\sigma}} & \leq & C \varepsilon^{\sigma} ,  \label{borne0} 
\\
|| u_{\varepsilon} ||_{\widetilde{B}^{-\sigma}} & \leq & C \varepsilon^{\sigma}  . \label{borne1}  \\
|| u_{\varepsilon} ||_{H^s} & \leq & C \varepsilon^{-s}   . \label{borne2}
\end{eqnarray}
\end{prop}

It follows from Theorem \ref{maintheorem} and  Proposition \ref{contreexemplevariete} with $ s = s_p $ (see  (\ref{pournotationsp})) and $ \sigma = \frac{n}{2} -s_p $ that (\ref{familyRn}) still holds for the family of functions (\ref{family}). In particular, we have the lower bound $ || u_{\varepsilon} ||_{H^{s_p}} \gtrsim \varepsilon^{-s_p} $. This shows that the definition of the norms $ || \cdot ||_{\widetilde{B}^{-\sigma}} $ is enough natural so that the estimates of Theorem \ref{maintheorem}  control accurately, as on $ \Ra^n $, certain fast oscillations which are badly estimated by pure Sobolev estimates. Note also that, in this counterexample, we did not require $u$ to be compactly supported so that Proposition \ref{contreexemplevariete} does not clearly follow from an elementary localization of the standard estimates on $ \Ra^n $.

\bigskip

\noindent {\it Proof.} We introduce first some notation (we refer to \cite{Bo1,Bo2} for more details). By possibly extending the function $w$ to the whole manifold, we can define a unitary mapping 
$$ L^2 \ni u \mapsto w(r)^{(1-n)/2} u \in \widehat{L^2} := L^2 ({\mathcal M}, w(r)^{(n-1)}d {\rm vol}_g) $$ and then set $ P = w(r)^{(1-n)/2}(- \Delta_g) w (r)^{(n-1)/2} $
which is selfadjoint on $ \widehat{L^2} $. 
We define $ v_{\epsilon}  $ to be the image of $ u_{\varepsilon} $ under this mapping. We next choose a fixed cutoff $ \varphi \in C^{\infty}({\mathcal M}) $ supported in $ (R,\infty) \times U $ which is equal to $1$ near the support of $ u $. We also choose local coordinates $ \kappa : U \rightarrow V \subset \Ra^{n-1} $ and define the operators $ J_{\kappa} : \widehat{L^2} \mapsto L^2 (\Ra^n) $ and $ J_{\kappa}^* : L^2 (\Ra^n) \rightarrow \widehat{L^2} $ by 
\begin{eqnarray*}
(J_{\kappa} v)(x_1,x^{\prime}) &= & \varphi (x_1, \kappa^{-1}(x^{\prime})) v (x, \kappa^{-1}(x^{\prime})) , \\
(J_{\kappa}^* f ) (r,m) & = & \varphi (r,m) f (r, \kappa (m)) .
\end{eqnarray*}
Considering the measure $ w(r)^{n-1} d {\rm vol_g} $ implies that these operators are bounded between $ L^p (\Ra^n) $ and $ \widehat{L^p} $ for all $p$.  In addition they satisfy the convenient relation $ J_{\kappa}^* J_{\kappa} v_{\varepsilon} = v_{\varepsilon} $.
 To prove (\ref{borne0}), we start by writing
\begin{eqnarray*}
 || u_{\varepsilon} ||_{H^{-\sigma}} & = & || (1 + P)^{-\sigma/2} v_{\epsilon} ||_{\widehat{L}^2}  \\
 & = & || (1 + P)^{-\sigma/2} J_{\kappa}^{-1} (1+D_1^2)^{\sigma/2} (1 + D_1^2)^{-\sigma/2} J_{\kappa} v_{\epsilon} ||_{\widehat{L}^2}
\end{eqnarray*} 
where, by  Theorem 1.5 of \cite{Bo2}, $ (1 + P)^{-\sigma/2} J_{\kappa}^{-1} (1+D_1^2)^{\sigma/2} $ is bounded from $ L^2 (\Ra^n) $ to $ \widehat{L^2} $. Here and below, we use the standard notation $ D_1 = i^{-1} \partial_{x_1} $. This implies that
\begin{eqnarray*}
 || u_{\varepsilon} ||_{H^{-\sigma}} & \leq & C \big| \big| (1+D_1^2)^{-\sigma/2} ( e^{i \frac{x_1}{\varepsilon}} \psi ) \big| \big|_{L^2 (\Ra)} , \\
 & \leq & C \big| \big| (1+|D_1 + \varepsilon^{-1} |)^{-\sigma/2}  \psi  \big| \big|_{L^2 (\Ra)} \\
 & \leq & C \varepsilon^{\sigma} ||  \psi||_{H^{\sigma}(\Ra)}  
\end{eqnarray*}
by using in the last estimate the Peetre inequality $ \scal{\xi_1 + \varepsilon^{-1}}^{-\sigma} \leq C \scal{\xi_1}^{\sigma} \varepsilon^{\sigma} $. To prove (\ref{borne1}), we next need to estimate $ || \theta (h A) u_{\varepsilon} ||_{L^{\infty}} $. Note that $ \theta (hA) = \widetilde{\theta} (h^2 \Delta_g) $ with $ \widetilde{\theta}(\lambda) = \theta (|\lambda|^{1/2}) $ which belongs to $ C_0^{\infty} (\Ra) $ so, by Theorem 1.5 of \cite{Bo2}, $ \theta (h A) $ is the sum of a pseudodifferential operator $ \mbox{Op}(\theta,h)  $ and a remainder which is bounded from $ H^{-\sigma} $ (for any $ \sigma > 0 $) to $ L^{\infty} $ so that
$$ || \theta (h A) u_{\varepsilon} ||_{L^{\infty}} \leq  C || \mbox{Op}(\theta,h) u_{\varepsilon} ||_{L^{\infty}} + C || u_{\varepsilon} ||_{H^{-\sigma}} , $$
with a constant $ C $ independent of $h$ (and $ u_{\varepsilon} $). To estimate $ || \mbox{Op}(\theta,h) u_{\varepsilon} ||_{L^{\infty}} $, it suffices to estimate the $ L^{\infty} (\Ra^n)$ norm of functions (of $ x= (x_1 , x^{\prime} ) $) of the form
$$ f_{\varepsilon} (x)  = e^{i \frac{x_1}{\varepsilon}} \int \! \! \int e^{i x_1 \xi_1 + i x^{\prime} \cdot \xi^{\prime}} a \big(x, h ( \xi_1 + \varepsilon^{-1} ) , h w (x_1) \xi^{\prime} \big) \widehat{\psi}(\xi_1) \widehat{\kappa_* \gamma}(\xi^{\prime})  d \xi_1  d \xi^{\prime} $$
where $ a$ is compactly supported and smooth in $ \xi $, and bounded in $x$. We have a first trivial bound
$$ || f_{\varepsilon} ||_{L^{\infty}(\Ra^n)} \leq C || \widehat{\psi} ||_{L^1} || \widehat{\kappa_* \gamma} ||_{L^1} \leq C \varepsilon^{\sigma} h^{-\sigma} $$
if $ h \leq \varepsilon $. When $ \varepsilon < h $, using the boundedness of $w$ and again the Peetre inequality, we observe that
\begin{eqnarray*}
 \big| a \big(x, h ( \xi_1 + \varepsilon^{-1} ) , h w (x_1) \xi^{\prime} \big) \big| & \leq & C (1 + |h \xi_1 + h \varepsilon^{-1}| + |h w (x_1) \xi^{\prime}|)^{-\sigma} \\
 & \leq & C (1+|h \xi_1| + |h w (x_1) \xi^{\prime}|)^{\sigma} (1 + |h \varepsilon^{-1}|)^{-\sigma} \\
 & \leq & C (1 + |\xi|)^{\sigma} h^{-\sigma} \varepsilon^{\sigma}
\end{eqnarray*} 
which implies that in this case we also have $ || f_{\varepsilon} ||_{L^{\infty}(\Ra^n)} \leq C h^{-\sigma} \varepsilon^{\sigma} $. All this implies precisely that $ || \theta (2^{-k}A) u_{\varepsilon}||_{L^{\infty}} \leq C \varepsilon^{\sigma} 2^{k \sigma} $, which is the meaning of (\ref{borne1}). Finally to prove (\ref{borne2}), we only need  to control the contribution of the pseudodifferential expansion of $ (1+P)^{s/2} v_{\varepsilon} $ in Theorem 1.5 of \cite{Bo2} which leads to estimate functions of the form 
$$  e^{i \frac{x_1}{\varepsilon}} \int \! \! \int e^{i x_1 \xi_1 + i x^{\prime} \cdot \xi^{\prime}} a_s \big(x,  \xi_1 + \varepsilon^{-1}  ,  w (x_1) \xi^{\prime} \big) \widehat{\psi}(\xi_1) \widehat{\kappa_* \gamma}(\xi^{\prime})  d \xi_1  d \xi^{\prime} $$
in $ L^2 (\Ra^n) $ with $a_s \in S^{s}$  a symbol of order $s$. In this integral, we write
$$ a_s \big(x,  \xi_1 + \varepsilon^{-1}  ,  w (x_1) \xi^{\prime} \big) = \left[ a_s \big(x,  \xi_1 + \varepsilon^{-1}  ,  w (x_1) \xi^{\prime} \big) \scal{\xi_1 + \varepsilon^{-1},\xi^{\prime}}^{-s} \right] \scal{\xi_1 + \varepsilon^{-1} , \xi^{\prime}}^s . $$
It is not hard to check that the bracket $ [ \cdots ] $ in the right hand side is a bounded family of symbols in $ S^0_0 $ as $ \varepsilon $ varies. By the Calder\'on-Vaillancourt Theorem, we thus get
\begin{eqnarray*}
 || u_{\varepsilon} ||_{H^s} &  \leq & C + C || \scal{\xi_1 + \varepsilon^{-1},\xi^{\prime}}^s \widehat{\psi}(\xi_1) \widehat{\kappa_* \gamma}(\xi^{\prime}) ||_{L^2 (\Ra^n_{\xi_1,\xi^{\prime}})} \\
 & \leq &  C \varepsilon^{-s} 
\end{eqnarray*} 
and the result follows. \finpreuve

\bigskip

\noindent {\bf Acknowledgements.} The first author thanks Michel Ledoux for interesting discussions and encouraging comments.

\end{document}